\documentclass[11pt]{amsart}
\usepackage{amssymb}
\usepackage{amsthm}
\usepackage{mathrsfs}
\usepackage[all]{xy}
\newtheorem{thm}{Theorem}[section]

\theoremstyle{definition}

\theoremstyle{remark}
\newtheorem{rem}[thm]{Remark}

\newcommand{\abs}[1]{\left\vert#1\right\vert}

\newcommand{\Cpx}{\mathbb{C}}

\newcommand{\sph}{\mathbb{S}}

\begin{document}

\title[Ricci solitons on Sasakian manifolds]{Ricci solitons on Sasakian manifolds}%
\author{Chenxu He and Meng Zhu}%

\address{14 E. Packer Ave \\
 Dept. of Math, Lehigh University \\
 Christmas-Saucon Hall \\
 Bethlehem, PA 18015}

\email{he.chenxu@lehigh.edu \\
mez206@lehigh.edu}


\begin{abstract}
We show that a Sasakian metric which also satisfies the gradient Ricci soliton equation is necessarily Einstein.
\end{abstract}

\maketitle

\section{introduction}
A Riemannian manifold $(M,g)$ is called a Ricci soliton if there exists a smooth vector field $X$ such that the Ricci tensor satisfies the following equation
\begin{equation}\label{eqnRicsoliton}
\mathrm{Ric} + \frac{1}{2}\mathscr{L}_X g = \lambda g
\end{equation}
for some constant $\lambda$ and $\mathscr{L}_X$ is the Lie derivative. It is called a gradient Ricci soliton if $X=\nabla f$ for some function $f$ on $M$. Ricci solitons appear as self-similar solutions to Hamilton's Ricci flow and often arise as limits of dilations of singularities in the Ricci flow, see \cite{Hamiltonformation}. They are also natural generalizations of Einstein metrics. Note that a soliton is called shrinking if $\lambda > 0$, steady if $\lambda=0$ and expanding if $\lambda <0$. It is well-known that when $\lambda \leq 0$ all compact solitons are necessarily Einstein.

The constructions of Ricci soliton metrics which are not Einstein are very limited. When the manifold is compact there are a few examples on K\"{a}hler manifolds due to N. Koiso,  H.-D. Cao, and X.-J. Wang and X. Zhu, see \cite{Koiso}, \cite{Caosoliton} and \cite{WangZhu}. In Koiso-Cao's construction, the manifolds are total spaces of $\Cpx \mathrm{P}^1$ bundles over a K\"{a}hler-Einstein manifold. Recently A. Dancer and M.-Y. Wang generalized Koiso-Cao's examples to allow the base as a product of manifolds using the framework of cohomogeneity one manifolds, see \cite{DancerWang}. Examples in low dimensions are even more rare. In dimension three or less all compact shrinking solitons are necessarily of positive constant curvature by R. Hamilton, T. Ivey and G. Perelman, see \cite{Hamiltonsurface}, \cite{Ivey3compact} and \cite{Perelman}. There are two compact non-Einstein examples in four dimension, $\Cpx \mathrm{P}^2 \sharp \overline{\Cpx \mathrm{P}^2}$ by Koiso and Cao, $\Cpx \mathrm{P}^2 \sharp 2\overline{\Cpx \mathrm{P}^2}$ by Wang-Zhu. Here $\overline{\Cpx \mathrm{P}^2}$ is a copy of $\Cpx \mathrm{P}^2$ with opposite orientation. All of these compact examples are K\"{a}hler manifolds. A natural question asked by H.-D. Cao is to find a non-product, non-Einstein, non-K\"{a}hler Ricci shrinking soliton, see \cite{CaoSurvey}.

As Ricci solitons are generalization of Einstein manifolds, one may look for new examples in the class where there are Einstein manifolds. A rich class of Einstein manifolds has been found in the class of Sasakian manifolds. Interesting examples of Sasakian-Einstein manifolds include exotic spheres, see \cite{BoyerGalicki} and references therein. Note that all examples of shrinking Ricci solitons are K\"{a}hler so far, thus of even dimensions, but Sasakian manifolds are of odd dimension. Moreover, a class of examples of Sasakian-Einstein metric is constructed on $\sph^1$-bundles over K\"{a}hler-Einstein manifolds. These considerations led Cao\cite{Caoprivate} to ask the following natural

Question: \emph{Does there exist a shrinking Ricci soliton on Sasakian manifolds which is
not Einstein?}

In the contrast to Einstein manifolds, we have the following non-existence result on gradient Ricci solitons.
\begin{thm}\label{thmmain}
Suppose $(M, g)$ is a Sasakian manifold and satisfies the gradient Ricci soliton equation,
then $f$ is a constant function and $(M,g)$ is an Einstein manifold.
\end{thm}

\begin{rem}Note that in our paper a Ricci soliton with Sasakian metric is different from the Sasaki-Ricci soliton in the context of transverse K\"{a}hler structure in a Sasakian manifold, see for example \cite{FutakiOnoWang}.
\end{rem}

\begin{rem}
Since all compact Ricci solitons are gradient ones from \cite{Perelman}, Theorem 1.1 implies that there is no compact non-Einstein Ricci soliton in Sasakian manifolds.
\end{rem}

\begin{rem}On the other hand, there exist non-gradient expanding Ricci solitons on non-compact Sasakian manifolds which are not Einstein. They were first discovered by P. Baird and L. Danielo in \cite{BairdDanielo}, and independently by J. Lott in \cite{Lott}. They are left invariant metrics on some solvable Lie groups and appear as type III singularities models of Ricci flow.  The one with the lowest dimension is the three dimensional Heisenberg group with a left invariant metric. One can easily check that the Reeb vector field is given by $F_3$ in Example 1.54 in \cite{Chowetc}.
\end{rem}

\begin{rem}
After the first version of our paper appeared in the arXiv, we learned that similar results have been obtained in $K$-contact manifolds by R. Sharma, see \cite{Sharma}.
\end{rem}

\smallskip

\textsc{Acknowledgment: } We thank Huai-Dong Cao for bringing this problem to our attention and useful discussions. We also want to thank Xiaofeng Sun for useful conservations, and Mircea Crasmareanu for informing us the reference \cite{Sharma}.

\medskip

\section{Proof of Theorem \ref{thmmain}}

In this section we prove Theorem \ref{thmmain}. The basic properties of Sasakian manifolds can be found in \cite{BoyerGalicki} or \cite[Chapter V]{YanoKon}.

\begin{proof}
Recall that a Sasakian manifold $(M^{2m+1}, g, \xi, \eta, \Phi)$ is a contact manifold on which there exists a unit Killing vector field $\xi$(called the Reeb vector field), $\eta$ is the dual one-form of $\xi$ and $\Phi$ is a $(1,1)$-tensor defined by $\Phi(Y) = \nabla_Y \xi$. The metric $(M, g)$ is called Sasakian if and only if
\begin{equation}\label{eqnSasakiR}
R(Y, \xi)Z = - g(Y, Z)\xi + g(Z, \xi) Y,
\end{equation}
for any vector fields $Y, Z \in TM$. Let $\mathcal{D} \subset TM$ be the distribution defined by $\eta(Y) = g(Y, \xi) = 0$. Then $\mathcal{D}$ is nowhere integrable as $\eta$ is a contact $1$-form. Since $\xi$ is a unit Killing vector field we have $\Phi(Y) = \nabla_Y \xi \in \mathcal{D}$ for any $Y \in TM$.

From the curvature condition (\ref{eqnSasakiR}) for any $Y \in \mathcal{D}$ and $Z \in TM$ we have
\begin{equation*}
R(Y, \xi, Z, Y) = - g(Y,Z)g(\xi,Y) + g(Z, \xi)g(Y,Y) = g(Z, \xi)\abs{Y}^2
\end{equation*}
and thus $\xi$ is an eigenvector field of $\mathrm{Ric}$ with
\begin{equation*}
\mathrm{Ric}(\xi) = 2m \xi.
\end{equation*}

In the following, we do the calculation on a general Ricci soliton. At the end, we use the fact that $X = \nabla f$ to draw the conclusion. Since $\xi$ is a Killing vector field, taking the Lie derivative to the soliton equation (\ref{eqnRicsoliton}) yields
\begin{equation*}
\mathscr{L}_{\xi}(\mathscr{L}_X g)= 0,
\end{equation*}

Thus, for any vector fields $Y, Z \in TM$ we have
\begin{eqnarray*}
0 &=& \left(\mathscr{L}_{\xi}(\mathscr{L}_X g)\right)(Y,Z) \\
& = & \mathscr{L}_{\xi}\left(\mathscr{L}_X g\right)(Y,Z) - \mathscr{L}_X g\left([\xi, Y],Z\right)-\mathscr{L}_X g\left(Y,[\xi,Z]\right) \\
& = & \nabla_{\xi}\left(\nabla_X g(Y,Z) - g\left([X,Y],Z\right) - g\left(Y, [X,Z]\right)\right) \\
& & - \nabla_X g\left([\xi, Y], Z\right) + g\left([X,[\xi, Y]], Z\right) + g\left([\xi, Y], [X, Z]\right) \\
& & - \nabla_X g\left( Y, [\xi, Z]\right)+ g\left([X, Y], [\xi, Z]\right) + g\left( Y, [X,[\xi, Z]]\right) \\
& = & \nabla_{\xi} \left(g\left( \nabla_Y X, Z\right) + g\left(\nabla_Z X, Y\right) \right) \\
& &  - g\left(\nabla_X [\xi, Y], Z\right) + g\left(\nabla_X [\xi, Y]- \nabla_{[\xi, Y]} X, Z\right) - g\left([\xi, Y], \nabla_Z X\right) \\
& & - g\left( Y, \nabla_X [\xi, Z]\right) - g\left(\nabla_Y X, [\xi, Z]\right) + g\left( Y, \nabla_X [\xi, Z] - \nabla_{[\xi, Z]} X\right) \\
& = & g\left(\nabla_{\xi}\nabla_Y X, Z\right) + g\left(\nabla_{\xi}\nabla_Z X, Y\right) - g\left(\nabla_{[\xi, Y ]} X, Z\right)\\
& & + g\left(\nabla_Y \xi, \nabla_Z X\right) - g\left(\nabla_{[\xi, Z ]}X, Y\right) + g\left(\nabla_Y X, \nabla_Z \xi\right), \\
& = &  R\left(\xi, Y, X, Z\right) +  g\left(\nabla_Y\nabla_{\xi}X, Z\right) +  R\left(\xi, Z, X, Y\right) +  g\left(\nabla_Z\nabla_{\xi}X, Y\right)\\
& & +  g\left(\nabla_Y \xi, \nabla_Z X\right) +  g\left(\nabla_Y X, \nabla_Z \xi\right).
\end{eqnarray*}
We let $Y$ be orthogonal to $\xi$ and $Z=\xi$ and then the above equation with the curvature condition (\ref{eqnSasakiR}) implies that
\begin{eqnarray}
0 \nonumber &=&  R\left(\xi, Y, X, \xi\right) +  g\left(\nabla_Y\nabla_{\xi} X, \xi\right) + g\left(\nabla_{\xi}\nabla_{\xi} X, Y\right) +  g\left(\nabla_Y \xi, \nabla_{\xi} X\right) +  g\left(\nabla_Y X, \nabla_{\xi}\xi\right)\\
&=& R\left(X, \xi, \xi, Y\right) + \nabla_{Y} g\left(\nabla_{\xi}X, \xi\right) + g\left(\nabla_{\xi}\nabla_{\xi} X, Y\right), \label{eqnR}
\end{eqnarray}
where the vector field $\nabla_{\xi}\xi$ vanishes as $\xi$ is a unit Killing vector field.

Note that from the soliton equation, for any $W \in TM$ we have
\begin{equation*}
2(\lambda - 2m)g(\xi, W) = 2(\lambda g(\xi, W) - \mathrm{Ric}(\xi, W)) = \left(\mathscr{L}_X g\right)(\xi, W) = g\left(\nabla_{\xi} X, W\right) + g\left(\nabla_W X, \xi \right).
\end{equation*}
Let $W = \xi$ then from the above equation we have
\begin{equation*}
2(\lambda - 2m) = 2g(\nabla_{\xi}X, \xi),
\end{equation*}
i.e,
\begin{equation*}
g(\nabla_{\xi}X, \xi) = \lambda - 2m,
\end{equation*}
or equivalently,
\begin{equation*}
\nabla_{\xi}g(X, \xi) = \lambda - 2m.
\end{equation*}
The middle term in the equation (\ref{eqnR}) vanishes and we have
\begin{eqnarray}
0 & = & R\left(X, \xi, \xi, Y\right) + g\left( \nabla_{\xi}\nabla_{\xi} X, Y\right) \notag \\
& = & g(X, Y) + g\left( \nabla_{\xi}\nabla_{\xi} X, Y\right), \label{eqnRtwoxi}
\end{eqnarray}
for any vector field $Y \perp \xi$.

In the case where $X = \nabla f$ for some smooth function, i.e., the soliton is a gradient one, we have $g(\nabla_{\xi}\nabla f, W) = (\lambda - 2m)g(\xi, W)$, i.e., $\nabla_{\xi}\nabla f = (\lambda - 2m)\xi$. So the second term in the equation (\ref{eqnRtwoxi}) vanishes as $\nabla_{\xi}\xi = 0$ and thus we have $g(\nabla f, Y) = 0$ for any $Y \perp \xi$, i.e., $\nabla f$ is parallel to $\xi$. Hence $\nabla f = 0$ as $\mathcal{D}$ is nowhere integrable, i.e., $f$ is a constant function.
\end{proof}



\medskip

\end{document}